\documentclass[12pt,a4]{article}
\usepackage{amsmath, amsthm, amssymb}
\usepackage[all,cmtip,pdf]{xy}
\usepackage{tikz-cd}
\usepackage{tikzsymbols}
\usepackage{amscd}
\usepackage{mathtools}
\usepackage[colorlinks=true,urlcolor=blue,citecolor=green,linkcolor=blue,bookmarks=true]{hyperref}

\newtheorem{theorem}{Theorem}[section]

\newtheorem{corollary}[theorem]{Corollary}

\newtheorem{conjecture}{Conjecture}

\newenvironment{conjecture*}{\begin{trivlist}\item[\hskip\labelsep{\bf Conjecture(*)}]}{%
\rm\end{trivlist}}

\newenvironment{proof*}{\begin{trivlist}\item[\hskip\labelsep{\bf Proof of Theorem \ref{t3.1}.}]}{%
\hfill$\square$\rm\end{trivlist}}
\newenvironment{proof**}{\begin{trivlist}\item[\hskip\labelsep{\bf Proof of Corollary \ref{c3.1}.}]}{%
\hfill$\square$\rm\end{trivlist}}
\newenvironment{ack}{\begin{trivlist}\item[\hskip\labelsep{\bf Acknowledgements.}]}{%
\hfill$\square$\rm\end{trivlist}}

\everymath{\displaystyle}

\begin{document}
\bibliographystyle{plainnat}
\pagestyle{plain} \pagenumbering{arabic}
\date{}
\newpage
\title{\textbf{\small{On The Matlis Duals of Local Cohomology Modules}}}
\author{\small{Gennady Lyubeznik and Tu\u{g}ba Y\i{Ä±}ld\i{Ä±}r\i{Ä±}m}\thanks{The first author gratefully acknowledges NSF support through grant DMS-1500264. The second author was supported by TÜBİTAK 2214/A Grant Program: 1059B141501072}}
\renewcommand\footnotemark{}
\maketitle

\begin{abstract}
Let $(R,\mathfrak{m})$ be a Noetherian regular local ring containing a field of characteristic $p>0$ and $I$ a nonzero ideal of $R$. In this short note, we prove that if $\operatorname{H}^i_I(R)\neq 0$, then $ \operatorname{Supp}_R(D(\operatorname{H}^i_{I}(R)))=\operatorname{Spec}(R)$.
\end{abstract}
 {\bf \emph{Keywords}:}\emph{ Local cohomology, Matlis duality, F-modules\\
 2000 Mathematics Subject Classification.  	13D45, 13H05.}
\section{Introduction}

\par
Let $(R,\mathfrak{m})$ be a Noetherian local commutative ring with unity, $I$ an ideal of $R$ and $E:=E_R(R/{\mathfrak{m}})$ an $R$-injective hull of the residue field $R/{\mathfrak{m}}$. Then for any $R$-module $M$, we denote by $\operatorname{H}^i_I(M)$ the $i$-th local cohomology module of $M$ supported in $I$ and by $\operatorname{D}(M):=\operatorname{Hom}_R(M,E)$ the Matlis dual of $M$.\\\\
Suppose now that $\operatorname{H}^i_I(R)=0$ for all $i\neq c$ and let \textbf{\em x}= $\{x_1,x_2,...,x_c\}$ be a regular sequence in $I$.  Hellus [\cite{Hel}, Corollary 1.1.4] proved that $I$ is a set theoretic complete intersection ideal defined by $x_i$ if and only if  $x_i$ form a $\operatorname{D}(\operatorname{H}^c_I(R))$-regular sequence. Motivated by this result, Hellus studied the associated primes of Matlis duals of the top local cohomology modules and conjectured the following equality:

\begin{equation*}
\operatorname{Ass}_R(D(\operatorname{H}^c_{(x_1,x_2,\cdots,x_c)}(R)))=\{ \mathfrak{p}\in Spec(R)\ |\  \operatorname{H}^c_{(x_1,x_2,\cdots,x_c)}(R/{\mathfrak{p}})\neq 0 \}
\end{equation*}
It has been shown that this conjecture holds true in many cases; see eg. \cite{Hel1}, \cite{Hel2}, \cite{Hel-Sch}, \cite{Hel-Stu}.\\ Furthermore, Hellus proved that the above conjecture is equivalent to the following condition [\cite{Hel}, Theorem 1.2.3]:
\begin{itemize}
\item If $(R, \mathfrak{m})$ is a Noetherian local domain, $c\geq1$ and $x_1,x_2,\cdots,x_c\in R$, then the implication \\
$\operatorname{H}^c_{(x_1,x_2,\cdots,x_c)}(R)\neq 0$ $\implies$ $0\in \operatorname{Ass}_R(D(\operatorname{H}^c_{(x_1,x_2,\cdots,x_c)}(R))$\\
holds.
\end{itemize}
We conjecture that if $R$ is regular, then the above implication holds for all non-zero ideals independently of the number of generators, i.e.
 \begin{conjecture} \label{conj1} Let $(R, \mathfrak{m})$ be a Noetherian regular local ring and $I$ be a non-zero ideal of $R$. If $\operatorname{H}^i_{I}(R)\neq 0$, then  $0\in \operatorname{Ass}_R(D(\operatorname{H}^i_{I}(R)))$.
\end{conjecture}
Note that Conjecture \ref{conj1} is not true for non-regular rings. For a concrete example of a Noetherian local ring $(A,\mathfrak{m})$ of dimension $>1$ such that $\operatorname{H}^1_{\mathfrak{m}}(A)=A/{\mathfrak{m}}$, hence  $0\notin \operatorname{Ass}_R(D(\operatorname{H}_{\mathfrak{m}}^{1}(A)))$, see [\cite{Goto-Ogawa}, Example 2.4]. The authors would like to thank M. Asgarzadeh for bringing this example to our attention.\\
We prove the following:
\begin{theorem}\label{t3.1}
Let $(R,\mathfrak{m})$ be a complete Noetherian regular local ring containing a field of characteristic $p>0$ and $\mathcal{M}$ be an $F$-finite module such that $0\notin \operatorname{Ass}(\mathcal{M})$. Then $0\in \operatorname{Ass}(D(\mathcal{M})).$
\end{theorem}
We would like to point out that $0\notin \operatorname{Ass}(\mathcal{M})$ is a necessary condition  of Theorem \ref{t3.1}.  Indeed, $R$ itself is an $F$-finite module and $0\in \operatorname{Ass}({R})$ but $0\notin \operatorname{Ass}(D({R}))=\operatorname{Ass}(E)=\{ \mathfrak{m} \}$.\\

As an immediate consequence of Theorem \ref{t3.1}, we obtained the main result of this paper which establishes Conjecture \ref{conj1}  in the equicharacteristic $p>0$ case:
\begin{corollary}\label{c3.1}
Let $(R,\mathfrak{m})$ be a Noetherian regular local ring containing a field of characteristic $p>0$ and $I$ a non-zero ideal of $R$. If $\operatorname{H}^i_{I}(R)\neq 0$, then $\operatorname{Supp}_R(D(\operatorname{H}^i_{I}(R)))=\operatorname{Spec}(R)$. 
\end{corollary}

\section{Preliminaries}
Throughout, $R$ is a commutative Noetherian regular ring containing a field of characteristic $p>0$.\\
Let $R^{'}$ be the additive group of $R$ regarded as an $R$- bi-module with the usual left action and with the right $R$- action defined by $r^{'}r=r^pr^{'}$ for all $r\in R$ and $r^{'}\in R^{'}$. The Frobenius functor
\begin{equation*}
F: R-mod \longrightarrow R-mod
\end{equation*}
of Peskine-Szpiro \cite{PSz} is defined by
\begin{center}
$F(M)=R^{'}\otimes_RM$
\end{center}
\begin{equation*}
\begin{tikzcd}
    F(M \arrow{r}{h}  &N)= (R^{'}\otimes_RM\arrow{r}{id\otimes_Rh}  &R^{'}\otimes_RN)
\end{tikzcd}
\end{equation*}
for all $R$-modules $M$ and all $R$-module homomorphisms $h$, where $F(M)$ acquires its $R$-module structure via the left $R$-module structure on $R^{'}$.\\
The iteration of a Frobenius functor on $R$ leads one to the iterated Frobenius functors $F^i(-)$ which are defined for all $i\geq 1$ recursively by $F^1(-)=F(-)$ and $F^{i+1}=F\circ F^i(-)$ for all $i\geq 1$.\\
Note that the Frobenius functor $F(-)$ is exact [\cite{Kunz}, Theorem 2.1]; $F(R)\cong R$ and for any ideal $I$ of $R$, $F(R/I)=R/{I^{[p]}}$, where $I^{[p]}$ is the ideal of $R$ generated by $p$-th powers of all elements of $I$ [\cite{PSz}, I.1.3d].\\
Note also that for any Artinian module $N$, $F(D(N))=D(F(N))$ [\cite{Lyu97}, Lemma 4.1] and so $R=F(R)=F(D(E))=D(F(E))$ implies   
 $F(E)=E$. Then it follows from Remark 1.0.(f) of \cite{Lyu97} that for any finitely generated $R$-module $M$, $F(D(M))=D(F(M))$.\\
Now, for an $R$-module $M$, define a Frobenius map $\psi_M : M\longrightarrow F(M)$ on $M$ by $\psi_M (m):=1\otimes m\in F(M)$ for all $m\in M$.
It is worth pointing out that if $\operatorname{ann}(m)=I\subseteq R$, then $\operatorname{ann}(\psi_M(m))=I^{p}$.\\ 
An $\mathcal{F}$-module $\mathcal{M}$ is an $R$- module equipped with $R$-module isomorphism $\theta : \mathcal{M}\longrightarrow F(\mathcal{M})$ which we call the structure morphism. \\
A generating morphism of an $\mathcal{F}$ module $\mathcal{M}$ is an $R$-module homomorphism $\beta : M\longrightarrow F(M)$, where $M$ is some $R$-module, such that $\mathcal{M}$ is the limit of the inductive system in the top row of the commutative diagram
\begin{equation*}
\begin{tikzcd}
   M \arrow{r}{\beta} \arrow{d}[swap]{\beta} & F(M) \arrow{d}{F(\beta)}\arrow{r}{F(\beta)}& F^2(M) \arrow{d}{F^2(\beta)}\arrow{r}{F^2(\beta)}& \cdots & \\
F(M) \arrow{r}{F(\beta)}& F^2(M) \arrow{r}{F^2(\beta)}&   F^3(M)  \arrow{r}{F^3(\beta)}& \cdots &
\end{tikzcd}
\end{equation*}
and $\theta : \mathcal{M}\longrightarrow F(\mathcal{M})$, the structure isomorphism of $\mathcal{M}$, is induced by the vertical arrows in this diagram.\\
If $\beta$ is an injective map, then the exactness of $F$ implies that all maps in the direct limit system are injective, so that $M$ injects into $\mathcal{M}$. In this case, we shall refer to $\beta$ as a root morphism of $\mathcal{M}$, and $M$ as a root of $\mathcal{M}$. If $\mathcal{M}$ is an $\mathcal{F}$-module possesing a root morphism $\beta: M\longrightarrow \mathcal{M}$ with $M$ finitely generated, then we say that $\mathcal{M}$ is $\mathcal{F}$-finite. In particular, $R$, with any $F$-module structure, is an $\mathcal{F}$-finite module.
\section{Proofs}
We begin this section by giving the proof of Theorem \ref{t3.1}:
\begin{proof*}Since $\mathcal{M}$ is an $F$-finite $R$ module, it follows from Proposition 2.3. of \cite{Lyu97} that there exists a root morphism $\beta:M\rightarrow F(M)$ such that
\begin{equation*}
\begin{tikzcd}\mathcal{M}=\varinjlim (M\arrow{r}{\beta}& F(M)\arrow{r}{F(\beta)}&F^2(M)\arrow{r}{ F^2(\beta)}&\cdots).
\end{tikzcd}
\end{equation*}
Then applying Matlis dual functor $\operatorname{D}(-)=\operatorname{Hom_R(-,E(R/{\mathfrak m}))}$ to $\mathcal{M}$, we obtain
\begin{equation*}
\begin{tikzcd} D(\mathcal{M})=\varprojlim (D(M)&\arrow{l}[swap]{D(\beta)} D(F(M))&\arrow{l}[swap]{D(F(\beta))}D(F^2(M))&\arrow{l}[swap]{ D(F^2(\beta))}\cdots).
\end{tikzcd}
\end{equation*}
But then since Frobenius functor commutes with $D(-)$, we can write $D(\mathcal{M})$ as
\begin{equation*}
\begin{tikzcd} D(\mathcal{M})=\varprojlim (N&\arrow{l}[swap]{\alpha} F(N)&\arrow{l}[swap]{F(\alpha)}F^2(N)&\arrow{l}[swap]{ F^2(\alpha)}\cdots),
\end{tikzcd}
\end{equation*}
where $N=\operatorname{D}(M)$ and $\alpha=\operatorname{D}(\beta)$.  Note that  since $\beta$ is injective and $F$ is exact,  $F^k(\alpha)$ is surjective for all $k\geq 0$. Note also that since $M$ is a finitely generated $R$-module, $N$ is Artinian and so are all $F^k(N)$, $k>0$.\\
On the other hand, since $0\notin \operatorname{Ass}(\mathcal{M})$, $I=\operatorname{Ann}(M)=\operatorname{Ann}(N)$ is a nonzero ideal of $R$. Then it follows that $\operatorname{Ann}(F(N))=I^{[p]}$  and so $\operatorname{Ker}(\alpha: F(N)\rightarrow N)\neq 0 $. \\
Now we claim that there exists an element $n^{'}=(n^{'}_0,n^{'}_1,\cdots,n^{'}_k,\cdots)\in \operatorname{D}(\mathcal{M})$ such that $n^{'}_k\in F^k(N)$ and $F^{k-1}(\alpha)(n^{'}_k)=n^{'}_{k-1}$ and with the property that  $\operatorname{ann}(n^{\prime}_k)\subseteq \mathfrak m^{k}$ for all $k\geq 4$.  To construct such an element, first let $b_1\in \operatorname{Soc}(\operatorname{Ker}(\alpha))\subseteq F(N)$, where $\operatorname{Soc}(\operatorname{Ker}(\alpha)):=\operatorname{Ann}_{\operatorname{Ker}(\alpha)}(\mathfrak{m})$ denotes the socle of $\operatorname{Ker}(\alpha)$ and  define $b_k$, for all $k\geq 2$, inductively as the image of $b_{k-1}$ under the Frobenius map (defined in the preceding section) on $F^{k-1}(N)$, that is $b_k:=\psi_{F^{k-1}(N)}(b_{k-1})=1\otimes b_{k-1}\in F^k(N)$. Then by induction on $k$ (considering that $\operatorname{ann}(b_1)=\mathfrak{m}$ and $\operatorname{ann}(x)=I$ implies $\operatorname{ann}(\psi(x))=I^{[p]}$), we have $\operatorname{ann}(b_k)=\mathfrak{m}^{[p^{k-1}]}$. On the other hand, since $b_1\in \operatorname{Ker}(\alpha):= \operatorname{Ker}(F^0(\alpha))$, an easy induction argument shows that $b_k\in \operatorname{Ker}(F^{k-1}(\alpha))$ for all $k\geq 0$. For if $b_{k-1}\in \operatorname{Ker}(F^{k-2}(\alpha))$, then $F^{k-1}(\alpha)(b_k)=F^{k-1}(\alpha)(1\otimes b_{k-1})=1\otimes F^{k-2}(\alpha)(b_{k-1})=0$.\\
Let now $n^{\prime}_0$ be an element of $N$ and, for every $1\leq k\leq 3$, choose $n^{\prime}_{k}\in F^{k}(N)$ such that $n^{\prime}_{k-1}=F^{k-1}(\alpha)(n^{\prime}_{k})$. For $k\geq 4$, define $n_k$ in such a way that $F^{k-1}(\alpha)(n_k)=n^{'}_{k-1}$. Then,  either $\operatorname{ann}(n_k)\subseteq \mathfrak m^{k}$ or  $\operatorname{ann}(n_k+b_k)\subseteq \mathfrak m^{k}$. Indeed, if $\operatorname{ann}(n_k+b_k)\nsubseteq \mathfrak{m}^{k}$,  there exists an element $y\in \mathfrak m\setminus \mathfrak{m}^{k}$ such that $y(n_k+b_k)=0$. Then it follows that $\operatorname{ann}(n_k)\subseteq \operatorname{ann}(yn_k)=\operatorname{ann}(yb_k)$. Since  $y\in \mathfrak m\setminus \mathfrak{m}^{k}$ and $\operatorname{ann}(b_k)=\mathfrak m^{[p^{k-1}]}$, we have $\operatorname{ann}(n_k)\subseteq \operatorname{ann}(yb_k)\subseteq  m^{p^{k-1}-k}  $. To prove the fact that $\operatorname{ann}(yb_k)\subseteq  m^{p^{k-1}-k}  $, suppose on the contrary that there exists an element $z\in \operatorname{ann}(yb_k)$ such that $z\notin \mathfrak m^{p^{k-1}-k}$. Then clearly,  $yz\in \operatorname{ann}b_k$. On the other hand as $R\cong \kappa[[X_1,...,X_n]]$, $\kappa\cong R/{\mathfrak{m}}$ a field of characteristic $p>0$,  and $y\notin \mathfrak{m}^k$ and $z\notin \mathfrak m^{p^{k-1}-k}$, we may write 
\begin{center}
\hspace{-1.7cm}$y=\underbrace{\Sigma^{k-1}_{i=1}\alpha_i \textbf{X}^i}_f+\underbrace{\Sigma^{ \infty}_{i=k}\alpha_i\textbf{X}^i}_{f^{\prime}}$\\
$z=\underbrace{\Sigma^{p^{k-1}-k-1}_{j=1}\beta_j \textbf{X}^j}_g+\underbrace{\Sigma^{\infty}_{j=p^{k-1}-k}\beta_j\textbf{X}^j}_{g^{\prime}}$
\end{center}
where $f$ and $g$ are non-zero polynomials over $\kappa[X_1,...,X_n]$, $\textbf{X}^l=X_1^{l_1}X_2^{l_2}\cdots X_n^{l_n}$ for any positive integer $l=l_1+l_2+...+l_n$. Then $yz=fg+ fg^{\prime}+gf^{\prime}+g^{\prime}f^{\prime}$. Since $\kappa[X_1,...,X_n]$ is an integral domain and $f$ and $g$ are non-zero, so is $fg$ and clearly $0\neq \operatorname{deg}(fg)\leq p^{k-1}-k-1+k-1=p^{k-1}-2$ from which it follows that $yz\notin \mathfrak{m}^{p^{k-1}}$. But then this contradicts the fact that $yz\in \operatorname{ann}(b_k)=\mathfrak{m}^{[p^{k-1}]}$. Hence $\operatorname{ann}(n_k)\subseteq \operatorname{ann}(yb_k)\subseteq  m^{p^{k-1}-k}  $, as desired. Now since $p^{k-1}-k\geq k$ for all $k\geq 4$, we have $\operatorname{ann}(n_k)\subseteq \mathfrak m^{p^{k-1}-k}\subseteq \mathfrak m^{k}$. 
 Then define
 \begin{displaymath}
   n^{'}_k = \left\{
     \begin{array}{lr}
       n_k,\ \ \ \  \ \ \ \ \ \ \ \ \  \textrm{if}\ \operatorname{ann}(n_k)\subseteq \mathfrak m^{k},\\
       n_k+b_k, \ \ \ \  \ \ \  \textrm{otherwise}.
     \end{array}
   \right.
\end{displaymath}

Clearly, $n^{'}=(n^{'}_0,n^{'}_1,\cdots,n^{'}_k,\cdots)\in \operatorname{D}(\mathcal{M})$ and $\operatorname{ann}(n^{\prime}_k)\subseteq \mathfrak m^{k}$ for all $k\geq 4$. This proves the claim.\\
Finally, $\operatorname{ann}(n^{'})=0$ for if $z\in \operatorname{ann}(n^{'})$, then $z\in \operatorname{ann}(n^{'}_k)$ for all $k\geq 0$ which implies that $z\in \bigcap_{n\in \mathbb{N}}\mathfrak m^n=\{0\}$. This completes the proof of Theorem \ref{t3.1}.
\end{proof*}
The proof of Corollary \ref{c3.1} is an immediate consequence of Theorem \ref{t3.1}:

\begin{proof**}
Without loss of generality, we may, and do, assume that $R$ is complete [\cite{Hel}, Remark 4.1.1]. Since $R$ is an $F$-finite module, so are its all local cohomology modules and since $0\notin \operatorname{Ass}_R(\operatorname{H}^i_{I}(R))$ for any nonzero ideal $I$ of $R$, the result follows from Theorem \ref{t3.1}.
\end{proof**}
\begin{ack} The results in this paper were obtained while the second author visited the School of Mathematics at University of Minnesota. The second author would like to thank TUBITAK for their support and also to all members of Mathematics Department of University of Minnesota for their hospitality. 
\end{ack}

\flushleft
\small{{Department of Mathematics,} 
{University of Minnesota,}
{127 Vincent Hall, 206 Church St.,}
Minneapolis, MN 55455, USA\\
\textbf{e-mail:} gennady@math.umn.edu}

\flushleft
\small{{Department of Mathematics,} 
{Istanbul Technical University,}
{Maslak, 34469, Istanbul, Turkey.}\\
\textbf{e-mail:} tugbayildirim@itu.edu.tr}

\end{document}